\documentclass{amsart}

\newtheorem{theorem}{Theorem}
\newtheorem{lemma}[theorem]{Lemma}
\newtheorem{proposition}[theorem]{Proposition}

\theoremstyle{definition}

\theoremstyle{remark}

\newcommand{\R}{{\mathbb R}}

\newcommand{\C}{{\mathbb C}}

\begin{document}
\title{On convexity of solutions of ordinary differential equations.}
\author[M.~Keller-Ressel]{Martin Keller-Ressel}
\address{Vienna University of Technology, Wiedner Hauptstrasse 8-10, A-1040 Vienna, Austria}
\email{mkeller@fam.tuwien.ac.at}
\author[E.~Mayerhofer]{Eberhard Mayerhofer}
\address{Vienna Institute of Finance, University of Vienna and Vienna University of Economics
and Business Administration, Heiligenst\"adterstrasse 46-48, 1190
Vienna, Austria} \email{eberhard.mayerhofer@vif.ac.at}
\author[A.~G.~Smirnov]{Alexander G. Smirnov}
\address{I.~E.~Tamm Theory Department, P.~N.~Lebedev
Physical Institute, Leninsky prospect 53, Moscow 119991, Russia}
\email{smirnov@lpi.ru}
\thanks{E.M. gratefully acknowledges financial support from WWTF (Vienna Science and Technology Fund).
The research of A.G.S. was supported by the grant LSS-1615.2008.2.}
\subjclass[2000]{65L05, 26B25, 06F20}
\date{}

\begin{abstract}
We prove a result on the convex dependence of solutions of ordinary
differential equations on an ordered finite-dimensional real vector
space with respect to the initial data.
\end{abstract}

\maketitle
\section{Introduction}
\label{s1}

Let $E$ be a finite-dimensional real vector space ordered by a closed proper cone\footnote{A set $C$ in a real vector space $E$ is called a cone if
$\lambda C\subset C$ for any $\lambda>0$. A cone $C$ is said to be proper if $C+C\subset C$ and $C\cap (-C)=\{0\}$. A cone $C$ induces a partial
order on $E$ if and only if it is proper.}~$C$.

Let $0<T\leq\infty$, $U\subset E$ be a non-empty open set, and
$f\colon [0,T)\times U\to E$ be a locally Lipschitz continuous map.
For any $x\in U$, the differential equation
\begin{equation}\label{a1}
\dot\psi(t) = f(t,\psi(t))
\end{equation}
has a unique maximally extended solution $\psi_f(\cdot,x)$ satisfying $\psi_f(0,x)=x$. This solution is defined on a semi-interval $[0,\theta_f(x))$,
where $0<\theta_f(x)\leq T$. For any $t\geq 0$, we set $\mathcal D_f(t)=\{x\in U: t<\theta_f(x)\}$.

Let $D\subset E$. A map $g\colon D\to E$ is called quasi-monotone
increasing~\cite{Volkmann} if the implication
\[
x\leq y,\,l(x)=l(y) \Longrightarrow\,\,l(g(x))\leq l(g(y))
\]
holds for all $x,y\in D$ and $l\in C^*$, where $C^*=\{l\in E^*:
l(x)\geq 0\mbox{ for any }x\in C\}$ is the dual cone of $C$ ($E^*$
is the dual space of $E$). A map $g\colon D\to E$ is called convex
if $D$ is convex and
\begin{equation}\label{conv}
g(\lambda x+(1-\lambda)y)\leq \lambda g(x)+(1-\lambda)g(y)
\end{equation}
for all $x,y\in D$ and $\lambda\in [0,1]$. A set $D\subset E$ is
said to be order regular if the relations $x\in D$ and $y\leq x$
imply that $y\in D$.

Our aim is to prove the next theorem.

\begin{theorem}\label{t1}
Let $U\subset E$ be a nonempty order-regular convex open set. Let
$0<T\leq \infty$ and $f\colon [0,T)\times U\to E$ be a continuous
map. If $f(t,\cdot)$ is quasi-monotone increasing and convex for all
$t\in [0,T)$, then $\mathcal D_f(t)$ is convex for any $t \in
[0,T)$, and $\psi_f(t,\cdot)$ is convex thereon.
\end{theorem}

In the formulation of Theorem~\ref{t1}, we do not require the local Lipschitz continuity of $f$ because the latter is ensured by continuity and
convexity (see Lemma~\ref{lLip} below). Note that the quasi-monotonicity of $f$ is a sufficient but not necessary condition for Theorem~\ref{t1} to
hold. For example, if $f(t,x)=f(x)$ is a linear map, then $\psi_f(t,x)$ is linear and hence convex in $x$, but $f$ may be not quasi-monotone
increasing in this case. On the other hand, at least in the autonomous case $f(t,x)=f(x)$, the convexity of $f$ is necessary to maintain the validity
of Theorem~\ref{t1}. Indeed, let $f$ be locally Lipschitz, $x,y\in U$ and $z=\lambda x+(1-\lambda)y$ with $0\leq \lambda\leq 1$. Suppose $\mathcal
D_f(t)$ is convex for any $t \in [0,T)$, and $\psi_f(t,\cdot)$ is convex thereon. Then we have
\[
\frac{\psi_f(t,z)-z}{t}\leq \lambda\frac{\psi_f(t,x)-x}{t}+(1-\lambda)\frac{\psi_f(t,y)-y}{t}
\]
for $t$ small enough. Passing to the limit $t\to 0$ in this inequality, we get $f(z)\leq \lambda f(x)+(1-\lambda)f(y)$, i.e., $f$ is convex.

The question of convex dependence of solutions of~(\ref{a1}) on initial data was first addressed in~\cite{Sarychev}, and then pursued
in~\cite{lakshawal96,Herzog2008}. In the last two papers, $E$ was assumed to be an ordered Banach space and it was shown (for differentiable $f$
in~\cite{lakshawal96} and for general locally Lipschitz continuous $f$ in~\cite{Herzog2008}) that $\psi_f(t,\cdot)$ is convex on any convex domain
contained in $\mathcal D_f(t)$ (in Appendix~A to this paper, we give a very simple proof of this result). Here, we strengthen this result in the
finite-dimensional case by proving the convexity of $\mathcal D_f(t)$. Moreover, keeping in mind possible applications (see, e.g., an example in
Section~\ref{s5}), we consider arbitrary open convex order-regular domains $U$ rather than the case $U=E$ studied in~\cite{lakshawal96, Herzog2008}.

The paper is organized as follows. In Section~\ref{s2}, we show that the conditions imposed on $f$ in Theorem~\ref{t1} ensure its local Lipschitz
continuity. In Section~\ref{s3}, we prove Theorem~\ref{t1} in the case, where $f$ is differentiable in the second variable. For this, we combine the
technique developed in~\cite{lakshawal96} with the well-known ``blow-up property'' of ODEs in finite dimensions: as $t\to \theta_f(x)$ for some $x\in
U$, the maximal solution $\psi_f(t,x)$ of~(\ref{a1}) must approach the boundary of the domain $[0,T)\times U$ on which $f$ is defined. In
Section~\ref{s4}, we get rid of the differentiability assumption and prove Theorem~\ref{t1} in the general case. Finally, in Section~\ref{s5}, we
illustrate Theorem~\ref{t1} by a concrete example of ODEs naturally arising in the theory of stochastic processes.

\section{Convexity and local Lipschitz continuity}
\label{s2}

Let $0<T\leq\infty$ and $\|\cdot\|$ be a norm on $E$. Let $U\subset
E$ be a non-empty open set. Recall that a map $f\colon [0,T)\times
U\to E$ is called locally Lipschitz if
\begin{equation}\label{LLL}
L_{t,K}(f) = \sup_{0\leq\tau\leq t,\,x_1,x_2\in K,\,x_1\ne x_2}
\frac{\|f(\tau,x_2)-f(\tau,x_1)\|}{\|x_2-x_1\|} <\infty
\end{equation}
for any compact set $K\subset U$ and any $t\in [0,T)$.

\begin{lemma} \label{lLip}
Let $f\colon [0,T)\times U\to E$ be a continuous map such that
$f(t,\cdot)$ is convex on $U$ for all $t\in [0,T)$. Then $f$ is
locally Lipschitz continuous.
\end{lemma}
\begin{proof}
Since $C$ is closed and $C\cap (-C)=\{0\}$, the set
$C\setminus\{0\}$ is contained in an open half-space of $E$. This
implies that the dual cone $C^*$ has a nonempty interior
(see,~e.g.,~\cite{Vladimirov1979}, Section~I.4.4, Lemma~1). Let
$l_1,\ldots,l_n\in C^*$ be a basis of $E^*$. Let the real-valued
functions $f_1,\ldots,f_n$ on $[0,T)\times U$ be defined by the
relations $f_j(t,x)=l_j(f(t,x))$. Clearly, $f_j$ are continuous on
$[0,T)\times U$ and $f_j(t,\cdot)$ are convex on $U$ for any $t\in
[0,T)$. Let $e_1,\ldots,e_n\in E$ be the dual basis of
$l_1,\ldots,l_n$: $l_j(e_k)=\delta_{jk}$. Then we have
\[
f(t,x)=\sum_{j=1}^n f_j(t,x) e_j.
\]
Hence, it remains to prove that $f_j$ are locally Lipschitz
continuous, i.e., satisfy~(\ref{LLL}) with $\|\cdot\|$ in the
numerator replaced with $|\cdot|$. Clearly, it suffices to
check~(\ref{LLL}) in the case $K=B_{x,r}$, where $B_{x,r}\subset U$
is a closed ball of radius $r>0$ centered at $x\in U$. Let $r'>r$ be
such that $B_{x,r'}\subset U$. By the continuity of $f_j$, there is
$m>0$ such that $|f_j(\tau,x)|\leq m$ for any $\tau\in [0,t]$ and
$x\in B_{x,r'}$. By~\cite[Corollary 2.2.12]{Zalinescu2002}, we have
\[
|f_j(\tau,x_2)-f_j(\tau,x_1)|\leq
\frac{2m}{r'}\frac{r'+r}{r'-r}\|x_2-x_1\|
\]
for any $x_1,x_2\in B_{x,r}$ and $\tau\in [0,t]$. The lemma is
proved.
\end{proof}

\section{The differentiable case}\label{s3}

In the rest of the paper, we assume that $T\in (0,\infty]$ is fixed
and set $I=[0,T)$.

Our consideration is essentially based on the next comparison result that is a particular case of a more general theorem proved by
Volkmann~\cite{Volkmann} in the setting of normed vector spaces.
\begin{lemma}\label{l2a}
Let $U\subset E$ be an open set. Let $f\colon I\times U\to E$ be a continuous locally Lipschitz map such that $f(t,\cdot)$ is quasi-monotone
increasing on $U$ for all $t\in I$. Let $0<t_0\leq T$ and $x,y\colon [0,t_0)\to U$ be differentiable maps such that $x(0)\leq y(0)$ and
\[
\dot x(t)-f(t,x(t))\leq \dot y(t)-f(t,y(t)), \quad \quad 0\leq
t<t_0.
\]
Then we have $x(t)\leq y(t)$ for all $t\in [0,t_0)$.
\end{lemma}

In fact, this comparison statement is essentially equivalent to
quasi-mono\-to\-ni\-city~\cite{Uhl1998}, but the above formulation
is  enough for our purposes. The next lemma is a simple
generalization of a well-known result for scalar-valued convex
functions.
\begin{lemma}\label{l3}
Let $U\subset E$ be an open convex set. A differentiable function
$g\colon U\to E$ is convex on $U$ if and only if
\begin{equation}\label{a4}
g(y)-g(x)\geq g'(x)(y-x), \quad x,y\in U.
\end{equation}
\end{lemma}
\begin{proof}
Let $h=y-x$ and $\lambda\in (0,1)$. If $g$ is convex on $U$, then
\[
g(x+\lambda h)=g((1-\lambda)x+\lambda y)\leq (1-\lambda)g(x)+\lambda
g(y).
\]
This implies that
\[
\frac{g(x+\lambda h)-g(x)}{\lambda}\leq g(y)-g(x).
\]
In view of the closedness of $C$, passing to the limit $\lambda \to
0$ yields~(\ref{a4}). Conversely, let (\ref{a4}) hold and $z=\lambda
x+(1-\lambda)y$. Then we have
\[
g(x)-g(z)\geq -(1-\lambda)g'(z) h,\quad g(y)-g(z)\geq \lambda g'(z)
h.
\]
Multiplying the left and right estimates by $\lambda$ and
$1-\lambda$ respectively and summing the resulting inequalities, we
obtain~(\ref{conv}). The lemma is proved.
\end{proof}

For differentiable functions, we have the following characterization
of quasi-monotonicity~\cite[Theorem~5]{Herzog2001}.
\begin{lemma}\label{l4}
Let $U\subset E$ be open and convex. A differentiable function $g\colon U\to E$ is quasi-monotone increasing on $U$ if and only if the linear map
$g'(x)\colon E\to E$ is quasi-monotone increasing for any $x\in U$.
\end{lemma}

Suppose $f\colon I\times U\to E$ is a continuous map such that $f(t,\cdot)$ is differentiable on $U$ for all $t\in I$ and the derivative
$f'(t,\cdot)$ is continuous on $I\times U$ (here and below, $f'(t,\cdot)$ denotes the derivative of the map $x\to f(x,t)$ with respect to $x$ for
fixed $t$). Then $f$ is locally Lipschitz, and we have
\begin{equation}\label{Lip}
L_{t,K}(f)=\sup_{0\leq \tau\leq t,\,x\in K} \|f'(\tau,x)\|
\end{equation}
for any $t\in I$, and for any compact convex set $K\subset U$ with a
nonempty interior. Given $x\in U$ and $0\leq t<\theta_f(x)$, we
define the linear map $B^x(t)\colon E\to E$ by setting
\begin{equation}\label{B}
B^x(t) = f'(t,\psi_f(t,x)).
\end{equation}
For $x\in U$ and $y\in E$, we denote by $w^x_y(t)$ the solution of
the initial value problem
\begin{equation}\label{a6}
\dot w^x_y(t) = B^x(t)w^x_y(t),\quad 0\leq t<\theta_f(x),\quad
w^x_y(0)=y.
\end{equation}
Clearly, $w^x_y$ is linear in $y$. For the norm of $w^x_y$, we have
the standard bound (see, e.g.,~\cite{Hartman}, Chapter IV,
Lemma~4.1)
\begin{equation}\label{w}
\|w^x_y(t)\|\leq \|y\|\exp\left(\int_0^t\|B^x(\tau)\|\,
d\tau\right),\quad 0\leq t<\theta_f(x).
\end{equation}

\begin{lemma}\label{l5}
Let $U\subset E$ be a convex open set and $f\colon I\times U\to E$
be a continuous map such that $f(t,\cdot)$ is differentiable on $U$
for all $t\in I$ and the derivative $f'(t,\cdot)$ is continuous on
$I\times U$. Suppose $f(t,\cdot)$ is convex and quasi-monotone
increasing on $U$ for all $t\in I$. For any $x,y\in U$, we have
\begin{equation}\label{a7}
w^x_{y-x}(t)\leq \psi_f(t,y)-\psi_f(t,x)\leq w^y_{y-x}(t),\quad
0\leq t< t_0,
\end{equation}
where $t_0=\min(\theta_f(x),\theta_f(y))$.
\end{lemma}
\begin{proof}
It suffices to prove the left inequality in~(\ref{a7}) because it
implies the right one after interchanging $x$ and $y$. Let
$s(t)=\psi_f(t,y)-\psi_f(t,x)$. By Lemma~\ref{l3}, we have
\[
\dot s(t)=f(t,\psi_f(t,y))-f(t,\psi_f(t,x))\geq B^x(t)s(t),\quad
0\leq t<t_0.
\]
By Lemma~\ref{l4}, the map $B^x(t)$ is quasi-monotone increasing for
any $t\in [0,t_0)$ and, therefore, the desired inequality follows
from~(\ref{a6}) and Lemma~\ref{l2a}. The lemma is proved.
\end{proof}

Since $E$ is finite-dimensional, the closed ordering cone $C$ is
normal. In terms of the partial order induced by $C$, this means
that there exists $\mu_C>0$ such that the implication
\begin{equation}\label{a8}
0\leq x\leq y \Longrightarrow\,\,\|x\|\leq\mu_C\|y\|
\end{equation}
holds for all $x,y\in E$.

If $f$ is continuously differentiable in the second variable,
Theorem~\ref{t1} follows from the next lemma.

\begin{lemma}\label{l6}
Let $U$ and $f$ be as in Lemma~$\ref{l5}$ and suppose in addition
that $U$ is order-regular. Let $x,y\in U$, $\lambda\in [0,1]$, and
$z=\lambda x + (1-\lambda)y$. Let
$t_0=\min(\theta_f(x),\theta_f(y))$. Then we have $\theta_f(z)\geq
t_0$ and
\begin{equation}\label{conv1}
\psi_f(t,z)\leq \lambda \psi_f(t,x) + (1-\lambda)\psi_f(t,y),\quad
0\leq t < t_0.
\end{equation}
Let $0\leq t<t_0$ and $K\subset U$ be a compact convex set with a
nonempty interior such that $\psi_f(\tau,x)$ and $\psi_f(\tau,y)$
lie in $K$ for all $\tau\in[0,t]$. Then
\begin{equation}\label{bound}
\|\psi_f(t,z)\|\leq R_K\left[1+\mu_C e^{L_{t,K}(f)\,t}\right],
\end{equation}
where $R_K=\sup_{\xi\in K}\|\xi\|$.
\end{lemma}
\begin{proof}
Let $\tau_0 = \min(\theta_f(x),\theta_f(y),\theta_f(z))$. Since
$z-x=(1-\lambda)(y-x)$ and $z-y=-\lambda(y-x)$, it follows from
Lemma~\ref{l5} that
\begin{align}
& (1-\lambda)w^x_{y-x}(t)\leq \psi_f(t,z)-\psi_f(t,x)\leq (1-\lambda)w^z_{y-x}(t),\nonumber\\
& -\lambda w^y_{y-x}(t)\leq \psi_f(t,z)-\psi_f(t,y)\leq -\lambda
w^z_{y-x}(t),\nonumber
\end{align}
for any $0\leq t<\tau_0$. Multiplying the first and second
inequalities by $\lambda$ and $1-\lambda$ respectively and adding
the results, we get
\begin{equation}\label{a9}
-\lambda(1-\lambda)v(t)\leq \psi_f(t,z) - u(t)\leq 0,\quad 0\leq
t<\tau_0,
\end{equation}
where $u,v\colon [0,t_0)\to E$ are given by
\begin{equation}\label{a10}
u(t) = \lambda\psi_f(t,x)+(1-\lambda)\psi_f(t,y),\quad v(t)=
w^y_{y-x}(t)-w^x_{y-x}(t).
\end{equation}
In view of~(\ref{a8}), it follows from~(\ref{a9}) that
\begin{multline}\label{a11}
\|\psi_f(t,z)\|\leq \|u(t)\|+\|\psi_f(t,z)-u(t)\|\leq \\ \leq
\|u(t)\|+\mu_C\lambda(1-\lambda)\|v(t)\|,\quad 0\leq t<\tau_0.
\end{multline}
Suppose that $\tau_0<t_0$. Then we obviously have
$\tau_0=\theta_f(z)$. Since both $u$ and $v$ are continuous on
$[0,t_0)$, it follows from~(\ref{a11}) that $\psi_f(t,z)$ is bounded
on $[0,\theta_f(z))$. This implies that we can choose a sequence
$t_k\uparrow \tau_0$ such that $\psi_f(t_k,z)$ converge to some
$x_0\in E$ as $k\to\infty$. By~(\ref{a9}), we have
$\psi_f(t_k,z)\leq u(t_k)$ for all $k$. As $C$ is closed, it follows
that $x_0\leq u(\tau_0)$. We hence have $x_0\in U$ because $U$ is
order-regular and $u(\tau_0)\in U$ by the convexity of $U$. On the
other hand, we cannot have $x_0\in U$ because $\psi_f(t,z)$ is a
maximal solution and must approach the boundary of $I\times U$ as
$t\to \theta_f(z)$ (see~\cite{Hartman}, Chapter~II, Theorem~3.1).
This contradiction shows that
\begin{equation}\label{t0}
\tau_0 = t_0.
\end{equation}
Combining this relation with~(\ref{a9}) and~(\ref{a10}), we
obtain~(\ref{conv1}). Let $t\in [0,t_0)$ and $K\subset U$ be a
convex compact set with a nonempty interior such that both
$\psi_f(\tau,x)$ and $\psi_f(\tau,y)$ lie in $K$ for any $\tau\in
[0,t]$. It follows from~(\ref{a10}), (\ref{w}), (\ref{B}),
and~(\ref{Lip}) that
\[
\|v(t)\|\leq 2\|y-x\|e^{L_{t,K}(f)\,t}\leq 4R_Ke^{L_{t,K}(f)\,t}.
\]
In view of~(\ref{t0}), inserting this estimate and the obvious inequalities $\|u(t)\|\leq R_K$ and $\lambda(1-\lambda)\leq 1/4$ in~(\ref{a11})
yields~(\ref{bound}). The lemma is proved.
\end{proof}

\section{Proof of Theorem~$\ref{t1}$}\label{s4}

To pass from continuously differentiable to arbitrary continuous functions, we shall need some results concerning the continuous dependence of
solutions of~(\ref{a1}) on the map $f$. Recall that Eq.~(\ref{a1}) possesses a maximal solution satisfying a given initial condition if the function
$f\colon I\times U\to E$ is continuous. Note however that such a solution may be not unique if $f$ is not locally Lipschitz continuous.

The next lemma easily follows from Theorem~3.2 in Chapter~II
of~\cite{Hartman}.

\begin{lemma}\label{l1}
Let $U\subset E$ be open. Let $f,f_1, f_2,\ldots$ be continuous maps
from $I\times U$ to $E$. Suppose $f$ is locally Lipschitz and $f_n$
converge to $f$ uniformly on all compact subsets of $I\times U$. Let
$\psi_n\in C^1([0,\theta_n),U)$ be maximal solutions of
\begin{equation}\label{a2}
\dot\psi_n(t) = f_n(t,\psi_n(t))
\end{equation}
such that $\psi_n(0)$ converge to some $u\in U$ as $n\to \infty$. Then we have
\begin{equation}\label{a3}
\theta_f(u) \leq \varliminf \theta_n.
\end{equation}
Let $0\leq a<\theta_f(u)$ and $n_0$ be such that $\theta_n>a$ for
$n> n_0$. Then the sequence $\psi_{n_0+k}(t)$, $k=1,2,\ldots$,
converges to $\psi_f(t,u)$ uniformly on $[0,a]$ as $k\to\infty$.
\end{lemma}

\begin{lemma}
\label{l2} Let $U\subset E$ be open. Let $f,f_1, f_2,\ldots$ be
continuous maps from $I\times U$ to $E$. Suppose $f$ is locally
Lipschitz and $f_n$ converge to $f$ uniformly on compact subsets of
$I\times U$. Let $0<a<T$ and $\psi_n\in C^1([0,a],U)$ be solutions
of~$(\ref{a2})$ such that $\psi_n(0)$ converge to some $u\in U$ as
$n\to \infty$. If for some compact set $K\subset U$, $\psi_n(t)\in
K$ for all $t\in [0,a]$, then $\theta_f(u)>a$, and we have
$\psi_n(t)\to \psi_f(t,u)$ and $\dot\psi_n(t)\to \dot\psi_f(t,u)$
uniformly on $[0,a]$.
\end{lemma}
\begin{proof}
Since $f_n$ are uniformly bounded on the compact set $Q=[0,a]\times K$, Eq.~(\ref{a2}) implies that $\dot \psi_n$ are uniformly bounded. Hence,
$\psi_n$ are uniformly equicontinuous. By the Arzel\`a-Ascoli theorem, it follows that the sequence $\psi_n$ is relatively compact in $C[0,a]$. Let
$\psi_{n_k}$ be a subsequence of $\psi_n$ uniformly converging to a function $\psi$. Obviously, $\psi(0)=u$ and $\psi(t)\in K$ for $t\in [0,a]$. Fix
$\varepsilon>0$. Because $f$ is uniformly continuous on $Q$, there exists a $\delta>0$ such that $\|f(t,x_1)-f(t,x_2)\|<\varepsilon/2$ for any
$(t,x_i)\in Q$ such that $\|x_2-x_1\|<\delta$. Let $k_0$ be such that $\|\psi_{n_k}(t)-\psi(t)\|<\delta$ and $\|f_{n_k}(t,x)-f(t,x)\|<\varepsilon/2$
for all $(t,x)\in Q$ and $k\geq k_0$. Then we have
\begin{multline}
\|f_{n_k}(t,\psi_{n_k}(t))-f(t,\psi(t))\|\leq\\ \leq
\|f_{n_k}(t,\psi_{n_k}(t))-f(t,\psi_{n_k}(t))\|+\|f(t,\psi_{n_k}(t))-f(t,\psi(t))\|
<\varepsilon,\quad t\in [0,a], \nonumber
\end{multline}
for any $k\geq k_0$, and in view of (\ref{a2}), the sequence $\dot
\psi_{n_k}(t)$ converges to $f(t,\psi(t))$ uniformly on $[0,a]$. On
the other hand, the uniform convergence of $\dot \psi_{n_k}$ implies
that $\psi$ is continuously differentiable and $\dot\psi$ is the
limit of $\dot \psi_{n_k}$. This means that $\psi$
satisfies~(\ref{a1}). Since $f$ is locally Lipschitz, this implies
that $\psi$ is the restriction of $\psi_f(\cdot,u)$ to $[0,a]$ and,
therefore, $\theta_f(u)>a$. We thus see that all uniformly
converging subsequences of $\psi_n$ have the same limit. As the
sequence $\psi_n$ is relatively compact, we conclude that
$\psi_n(t)\to\psi_f(t,u)$ uniformly on $[0,a]$. Replacing
$\psi_{n_k}$ with $\psi_n$ in the above proof, we obtain the uniform
convergence of $\dot \psi_n$. The lemma is proved.
\end{proof}

Proof of Theorem~$\ref{t1}$:

\begin{proof} For $\kappa>0$, we set
$U(\kappa)=\{\xi\in U: B_{\xi,\kappa}\subset U\}$, where $B_{\xi,\kappa}$ is the closed ball of radius $\kappa$ centered at $\xi$. Clearly, the set
$U(\kappa)$ is open, convex, and order-regular for any $\kappa>0$. Let $t\in I$, $x,y\in \mathcal D_f(t)$ and $z=\lambda x+(1-\lambda) y$ for some
$\lambda\in [0,1]$. We have to show that $\theta_f(z)>t$ and inequality~(\ref{conv1}) holds. Let $S\subset U$ be a convex compact set whose interior
contains $\psi_f(\tau,x)$ and $\psi_f(\tau,y)$ for all $\tau\in [0,t]$. Choose $\kappa>0$ such that $S\subset U(\kappa)$.

Let $\rho$ be a nonnegative smooth function on $E$ such that $\rho(\xi)=0$ for $\|\xi\|>1$ and $\int_{E} \rho(\xi)\,d\xi =1$. For any positive
$\varepsilon\leq \kappa$, we define the map $f_\varepsilon\colon I\times U(\kappa)\to E$ by setting
\[
f_\varepsilon(\tau,\xi) = \int_E f(\tau,\xi-\varepsilon
\eta)\rho(\eta)\,d\eta.
\]
Let $\phi$ denote the restriction of $f$ to $I\times U(\kappa)$.
Clearly, $f_\varepsilon$ are smooth in the second variable and
converge to $\phi$ uniformly on compact subsets of $I\times
U(\kappa)$ as $\varepsilon\to 0$. It is straightforward to check
that $f_\varepsilon$ are convex quasi-monotone increasing maps on
$U(\kappa)$ such that
\begin{equation}\label{a12}
L_{t,S}(f_\varepsilon)\leq L_{t, S_{\kappa}}(f),
\end{equation}
where $S_\kappa$ is the closed $\kappa$-neighborhood of $S$. Our
choice of $\kappa$ ensures that
$t<\min(\theta_{\phi}(x),\theta_{\phi}(y))$. Let
$t_\varepsilon=\min(\theta_{f_\varepsilon}(x),\theta_{f_\varepsilon}(y))$.
By Lemma~\ref{l1}, there exists $0<\varepsilon_0\leq \kappa$ such
that $t_\varepsilon>t$ for any $0<\varepsilon\leq \varepsilon_0$ and
$\psi_{f_\varepsilon}(\cdot,x)\to \psi_f(\cdot,x)$ and
$\psi_{f_\varepsilon}(\cdot,y)\to \psi_f(\cdot,y)$ uniformly on
$[0,t]$ as $\varepsilon_0\geq \varepsilon\to 0$. Decreasing
$\varepsilon_0$ if necessary, we can ensure that
$\psi_{f_\varepsilon}(\tau,x)$ and $\psi_{f_\varepsilon}(\tau,y)$
lie in $S$ for all $\tau\in[0,t]$ and $\varepsilon\in
(0,\varepsilon_0]$. It follows from Lemma~\ref{l6} that
$\theta_{f_\varepsilon}(z)\geq t_\varepsilon>t$ and
\begin{align}
&\psi_{f_\varepsilon}(\tau,z) \leq \lambda
\psi_{f_\varepsilon}(\tau,x)+(1-\lambda)\psi_{f_\varepsilon}(\tau,y),\label{convex eq}\\
&\|\psi_{f_\varepsilon}(\tau,z)\|\leq R_S\left[1+\mu_C e^{t
L_{t,S_\kappa}(f)}\right]\label{a13}
\end{align}
for any $0\leq \tau\leq t$ and $0<\varepsilon\leq\varepsilon_0$. Let
$r>0$ and $K=(S-C)\cap \{\xi\in E: \|\xi\|\leq r\}$. Since $S$ is
compact and $C$ is closed, $S-C$ is closed and, therefore, $K$ is
compact. The order-regularity of $U(\kappa)$ implies that $K\subset
U(\kappa)$. By~(\ref{convex eq}) and~(\ref{a13}), we have
$\psi_{f_\varepsilon}(\tau,z)\in K$ for all $0\leq \tau\leq t$ if
$r$ is large enough. It follows from Lemma~\ref{l2} that
$\theta_\phi(z)>t$ and $\psi_{f_\varepsilon}(\cdot,z)\to
\psi_{\phi}(\cdot,z)$ uniformly on $[0,t]$. Obviously, $\mathcal
D_\phi\subset \mathcal D_f$ and $\psi_\phi$ is the restriction of
$\psi_f$ to $\mathcal D_\phi$. Hence $\theta_f(z)\geq
\theta_\phi(z)>t$ and passing to the limit $\varepsilon\to 0$ in
inequality~(\ref{convex eq}) for $\tau=t$ yields~(\ref{conv1}). The
theorem is proved.
\end{proof}

\section{Example}\label{s5}
As an illustration, we give an example of a system of ODEs that
arises naturally in the theory of stochastic processes and satisfies
all conditions of Theorem~\ref{t1}. We consider a so-called
\emph{affine process} evolving on the state space $C:=\R_{\ge 0}^d$
(see \cite{DFS2003}). Such a process $X = (X_t)_{t \geq 0}$, can be
regarded as a multi-type extension of the singe-type continuously
branching process of~\cite{Lamperti1967}, which arises as a
continuous-time limit of a classical Galton-Watson branching
process. $X$ is defined as a stochastically continuous,
time-homogeneous Markov process starting at $X_0\in C$, with the
property that the moment generating function is of the form
\begin{equation}\label{Eq:affine_def}
\mathbb{E}\left[e^{x\cdot X_t}\right] = e^{\psi(t,x) \cdot X_0}
\end{equation}
for all $(t,x) \in \R_{\ge 0}\times \R^d$, and where $\psi: \R_{\ge
0} \times \R^d \to \R^d \cup \{\infty\}$.\footnote{We set $\psi(t,x)
= \infty$, whenever the left side of \eqref{Eq:affine_def} is
infinite. Note that for $(t,x) \in \R_{\ge 0} \times (-\infty,0]^d$
it is always guaranteed that $\psi(t,x)$ is finite.} We assume that
the time-derivative of $\psi(t,x)$ at $t = 0$,
\[
f(x) := \left.\frac{\partial}{\partial t}\psi(t,x)\right|_{t = 0}
\]
exists and is a continuous function on the set $U = \{x \in \R^d: f(x) < \infty\}$. In this case the map $\psi(t,x)$ satisfies the following
differential equation:
\begin{equation}\label{Eq:Riccati}
\frac{\partial}{\partial t}\psi(t,x) = f(\psi(t,x)), \qquad
\psi(0,x) = x.
\end{equation}
Moreover, the components of the map $f(x)$ are of so-called Levy-Khintchine type (cf. \cite[Theorem~8.1]{Sato1999}):
\[
f_i(x) = \frac{\alpha_i}{2} x_i^2  + x \cdot \beta^i  - c_i +
\int_{C\setminus\{0\}}{\left(e^{x \cdot \xi} - 1 - x \cdot \xi
\mathbf{I}_{|\xi| \le 1}\right)\,\mu_i(d\xi)},
\]
with $\mathbf{I}$, the indicator function, where, for all $i \in
\{1, \dotsc, d\}$,
\begin{itemize}
\item $\alpha_i \in \R_{\ge 0}$;

\item $\beta^i \in \R^d$ with $\beta^i_k - \int_{|\xi| \le 1}{\xi_k\,\mu_i(d\xi)} \ge 0$ for all $k \neq i$;

\item $c_i \in \R_{\ge 0}$;

\item $\mu_i(d\xi)$ are Borel measures on $C\setminus\{0\}$ assigning finite mass to the set $\{\xi \in C: |\xi| > 1\}$ and satisfying the
integrability condition
\[
\int_{\xi\in C,0<|\xi| \le 1}{\left(\sum_{k \neq i}|\xi_k| +
|\xi_i|^2 \right)\,\mu_i(d\xi)} < \infty
\]
on its complement.
\end{itemize}
The above conditions are both necessary and sufficient for the
existence of $X$ and referred to as \emph{admissibility} conditions
(see \cite{DFS2003}).

In the following we consider the ordering on $\mathbb R^d$ induced
by the cone $\mathbb R_{\geq 0}^d$.
\begin{proposition}
The domain $U$ is convex and order-regular and the map $f(x)$ is convex and quasi-monotone increasing thereon.
\end{proposition}
\begin{proof}
We make use of the following representations of $f_i(x)$:
\begin{equation}\label{R_rep}
f_i(x) = \log \int_{\R^d}{e^{x \cdot \xi}\,p_i(d\xi)} = f^\dagger_i(x) +\int_{C\setminus\{0\},\,|\xi|
> 1}{\left(e^{x
\cdot \xi} - 1\right)\,\mu_i(d\xi)},
\end{equation}
where $p_i(d\xi)$ is an infinitely divisible, substochastic measure
on $\R^d$, and $f^\dagger_i(x)$ is a function on $\R^d$, that can be
extended to an entire function on $\C^d$. The representation as
$\log \int_{\R^d}{e^{x \cdot \xi}\,p_i(d\xi)}$ is an immediate
consequence of the Levy-Khintchine formula, and its analytic
extension to exponential moments \cite[Theorem~8.1,
Theorem~25.17]{Sato1999}. The second representation of $f_i(x)$
follows directly from \cite[Lemma~25.6]{Sato1999}. To show that
$f_i(x)$ is convex, apply H\"older's inequality:
\begin{multline*}
f_i\left(\lambda x + (1 - \lambda) y\right) = \log
\int_{\R^d}{e^{\lambda x \cdot \xi} e^{(1 - \lambda) y \cdot
\xi}\,p_i(d\xi)} \le \\ \le \lambda \log \int_{\R^d}{e^{x \cdot
\xi}\,p_i(d\xi)} + (1 - \lambda) \log \int_{\R^d}{e^{y \cdot
\xi}\,p_i(d\xi)} = \lambda f_i(x) + (1 - \lambda) f_i(y)
\end{multline*}
for all $x,y \in \R^d$ and $\lambda \in (0,1)$. We show next that
the domain $U$ is order-regular. Assume that $x \in U$, i.e. $f_i(x)
< \infty$ for all $i$, and let $y \le x$. Using the second
representation in \eqref{R_rep} it is clear that $f_i^\dagger(y) <
\infty$. But also the integral with respect to $\mu_i(d\xi)$ is
finite, since the integrand is dominated by $(e^{x \cdot \xi} -
1)\mathbf{1}_{|\xi| \ge 1}$, whose integral is finite by assumption.
We conclude that $f_i(y) < \infty$, and thus that $y \in U$, i.e.,
$U$ is order-regular. Finally we show that $f(x)$ is also
quasi-monotone increasing. Assume that $y \le x$ with $y_i = x_i$
for some $i \in \{1, \dotsc, d\}$. It follows that
\begin{multline}
f_i(x) - f_i(y) = \sum_{k \neq i}{(x_k - y_k)\cdot \left(\beta^i_k -
\int_{\xi\in C,\,0<|\xi| \le 1}{\xi_k\,\mu_i(d\xi)}\right)} +\\+
\int_{C}{\left(e^{x \cdot \xi} - e^{y \cdot
\xi}\right)\,\mu_i(d\xi)} \ge 0, \nonumber
\end{multline}
where we have made use of the admissibility conditions given above.
\end{proof}

\section*{Appendix~A}
In this section we give a very simple proof of the convexity result~\cite{Herzog2008} for ODEs in ordered normed spaces. Let $E$ be a real normed
space (not necessarily finite-dimensional) ordered by a proper closed cone $C$. As shown in~\cite{Volkmann}, Lemma~\ref{l2a} holds for $E$ if one of
the following conditions is satisfied:
\begin{enumerate}
\item $C$ has a non-empty interior,

\item $E$ is complete,

\item $C$ is a distance set (i.e., for every $x\in E$, there is $y\in C$ such that $\|x-y\|$ is equal to the distance from $x$ to $C$).
\end{enumerate}

As above, let $T\in (0,\infty]$ and $I=[0,T)$. Theorem~1 in~\cite{Herzog2008} follows immediately from the next result.

\begin{theorem}
Let $E$ be an ordered normed space such that one of the above conditions is satisfied. Let $U\subset E$ be an open convex set and $f\colon I\times
U\to E$ be a continuous locally Lipschitz map such that $f(t,\cdot)$ is quasi-monotone increasing and convex on $U$ for all $t\in I$. Let $0<t_0\leq
T$ and $x_1,x_2,x_3\colon [0,t_0)\to U$ be differentiable maps such that
\[
\dot x_i(t)=f(t,x_i(t)),\quad i=1,2,3,
\]
and $x_3(0)=\lambda x_1(0)+(1-\lambda)x_2(0)$ for some $\lambda\in [0,1]$. Then $x_3(t)\leq \lambda x_1(t)+(1-\lambda)x_2(t)$ for all $t<t_0$.
\end{theorem}
\begin{proof}
Set $z(t)=\lambda x_1(t)+(1-\lambda)x_2(t)$ for $t<t_0$. By the convexity of $f$,
\begin{align*}
&\dot z(t)-f(t,z(t))= \lambda \dot x_1(t)+(1-\lambda)\dot x_2(t)-f(t,\lambda x_1(t)+(1-\lambda)x_2(t))\geq\\
&\geq \lambda(\dot x_1(t)-f(t,x_1(t))+(1-\lambda)(\dot x_2(t)-f(t,x_2(t)))=0=\dot x_3(t)-f(t,x_3(t))
\end{align*}
for all $t<t_0$. Since $z(0)=x_3(0)$, the above-mentioned analogue of Lemma~\ref{l2a} for normed spaces implies that $z(t)\geq x_3(t)$. The theorem
is proved.
\end{proof}


\begin{thebibliography}{10}
\expandafter\ifx\csname url\endcsname\relax
  \def\url#1{\texttt{#1}}\fi
\expandafter\ifx\csname urlprefix\endcsname\relax\def\urlprefix{URL }\fi

\bibitem{DFS2003} D.~Duffie, D.~Filipovi{\'c}, W.~Schachermayer, Affine processes and
  applications in finance, Ann. Appl. Probab. 13~(3) (2003) 984--1053.

\bibitem{Hartman} P.~Hartman, Ordinary differential equations, John Wiley \& Sons Inc., New York,
  1964.

\bibitem{Herzog2001} G.~Herzog, Quasimonotonicity, Nonlinear Anal. 47~(4) (2001) 2213--2224.

\bibitem{Herzog2008} G.~Herzog, R.~Lemmert, A note on convex dependence of solutions of {IVP}s
  relative to initial values, Nonlinear Anal. 68~(12) (2008) 3841--3844.

\bibitem{lakshawal96} V.~Lakshmikantham, N.~Shahzad, W.~Walter, Convex dependence of solutions of
  differential equations in a {B}anach space relative to initial data,
  Nonlinear Anal. 27~(12) (1996) 1351--1354.

\bibitem{Lamperti1967} J.~Lamperti, Continuous state branching processes, Bull. Amer. Math. Soc. 73
  (1967) 382--386.

\bibitem{Sarychev} A.~V. Sarychev, On equation {$x\sp {(n+1)}=f(t,x,\dot x,\cdots,x\sp {(n)})$}
  with convex quasi-monotone right-hand side, Nonlinear Anal. 27~(7) (1996)
  785--792.

\bibitem{Sato1999} K.~Sato, L\'evy processes and infinitely divisible distributions, vol.~68 of
  Cambridge Studies in Advanced Mathematics, Cambridge University Press,
  Cambridge, 1999, translated from the 1990 Japanese original, Revised by the
  author.

\bibitem{Uhl1998} R.~Uhl, Ordinary differential inequalities and quasimonotonicity in ordered
  topological vector spaces, Proc. Amer. Math. Soc. 126~(7) (1998) 1999--2003.

\bibitem{Vladimirov1979} V.~S. Vladimirov, Generalized functions in mathematical physics, ``Mir'',
  Moscow, 1979, translated from the second Russian edition by G.
  Yankovski{\u\i}.

\bibitem{Volkmann} P.~Volkmann, \"{U}ber die {I}nvarianz konvexer {M}engen und
  {D}ifferentialungleichungen in einem normierten {R}aume, Math. Ann. 203
  (1973) 201--210.

\bibitem{Zalinescu2002} C.~Z{\u{a}}linescu, Convex analysis in general vector spaces, World Scientific
  Publishing Co. Inc., River Edge, NJ, 2002.

\end{thebibliography}
\end{document}